\documentclass[11pt,a4paper,twoside]{amsart}

\usepackage{amssymb}
\usepackage{graphics}
\usepackage{epic,eepic}

\newtheorem{theo}{\bf Theorem}[section]
\newtheorem{prop}[theo]{\bf Proposition}
\newtheorem{lemma}[theo]{\bf Lemma}
\newtheorem{coro}[theo]{\bf Corollary}
\theoremstyle{remark}
\newtheorem{remark}[theo]{\bf Remark}

\newcommand{\RB}{{\rm RB}}
\newcommand{\crs}{\circlearrowleft}

\newcommand{\id}{{\rm id}}
\newcommand{\kors}{\#}
\DeclareMathOperator{\Inv}{Inv}
\DeclareMathOperator{\inv}{inv}
\DeclareMathOperator{\nega}{neg}
\DeclareMathOperator{\poin}{Poin}
\providecommand{\abs}[1]{\lvert#1\rvert}
\newcommand{\fs}{{\mathfrak{S}}}
\newcommand{\DEF}{\buildrel {\rm def} \over =}
\newcommand{\mat}[4]{{\bigl(\begin{smallmatrix}{#1} & {#2} \\
      {#3} & {#4}\end{smallmatrix}\bigr)}}

\begin{document}

\title{Bruhat intervals as rooks on skew Ferrers boards}
\author{Jonas Sj{\"o}strand}
\address{Department of Mathematics, Royal Institute of Technology \\
  SE-100 44 Stockholm, Sweden}
\email{jonass@kth.se}
\date{25 January 2006}

\begin{abstract}
We characterise the permutations $\pi$ such that the elements in the
closed lower Bruhat interval $[\id,\pi]$ of the symmetric group
correspond to non-taking
rook configurations on a skew Ferrers board.
It turns out that these are exactly the permutations $\pi$
such that $[\id,\pi]$ corresponds to a flag manifold defined by
inclusions, studied by Gasharov and Reiner.

Our characterisation connects the Poincar\'e polynomials
(rank-generating function)
of Bruhat intervals with $q$-rook polynomials, and
we are able to compute the Poincar\'e polynomial
of some particularly interesting intervals in the
finite Weyl groups $A_n$ and $B_n$. The expressions
involve $q$-Stirling numbers of the second kind.

As a by-product of our method, we present
a new Stirling number identity connected to both
Bruhat intervals and the poly-Bernoulli numbers defined
by Kaneko.
\end{abstract}

\maketitle

\section{Introduction}
Since its introduction in the 1930s
the Bruhat order on Coxeter
groups has attracted mathematicians from many areas.
Geometrically it describes the containment ordering of Schubert
varieties in flag manifolds and other homogeneous spaces.
Algebraically it is intimately related to the representation theory
of Lie groups. Combinatorially the Bruhat order is essentially the
subword order on reduced words in the alphabet of generators of a
Coxeter group.

The interval structure of the Bruhat order is geometrically
very important and has been studied a lot in the literature.
From a combinatorial point of view, as soon as there is a
(graded) poset, the following three questions naturally arise about
its intervals $[u,w]$ (and they will probably arise
in the following order):
\begin{enumerate}
\item
What is the rank-generating function (or Poincar\'e polynomial)
$\poin_{[u,w]}(q)=\sum_{v\in[u,w]}q^{\ell(v)}$?
\item
What is the M\"obius function $\mu(u,w)$?
\item
What can be said about the topology of the order complex of $(u,w)$?
\end{enumerate}
The third question was answered by Bj\"orner and Wachs ~\cite{bjornerwachs}
in 1982: The order complex of an open interval
$(u,w)$ is homeomorphic to the sphere $\mathbb{S}^{\ell(u,w)-2}$.
The second question was answered already by Verma~\cite{verma}
in 1971: $\mu(u,w)=(-1)^{\ell(u,w)}$. However, the
first question is still a very open problem!

For the whole poset $\poin(q)$ was computed
by Steinberg~\cite{steinberg}, Chevalley~\cite{chevalley},
and Solomon~\cite{solomon}. Really small intervals (of length
$\le7$ in $A_n$ and $\le5$ in $B_n$ and $D_n$) were completely classified
by Hultman~\cite{hultman} and Incitti~\cite{incitti}. Lower intervals
of 312-avoiding permutations in $A_n$ were classified by
Develin~\cite{develin} (though he did not compute
their Poincar\'e polynomials), and for a general
lower interval $[\id,w]$ in a crystallographic Coxeter group, Bj\"orner and
Ekedahl~\cite{bjornerekedahl} showed that the coefficients
of $\poin_{[\id,w]}(q)$ are partly increasing. Apart
from this, virtually nothing is known, not even for lower
intervals.

The aim of this paper is to start filling the hole and at least gain
some understanding of the rank-generating function of a family of
intervals in finite Weyl groups. To this end we present a connection
between the Poincar\'e polynomial $\poin(q)$
and rook polynomials, making it possible to compute $\poin(q)$
for various interesting intervals.
Our approach is partly a generalisation of the notion of
{\em partition varieties} introduced by Ding~\cite{ding} to
what may be called {\em skew partition varieties}.

The paper is composed as follows. In Section~\ref{sec:rook} we give a
short introduction to rook polynomials
before presenting our results in
Section~\ref{sec:results}. In Section~\ref{sec:connection}
we present the connection between rook polynomials and
Poincar\'e polynomials and prove our main theorem.
In sections~\ref{sec:typeAq} and~\ref{sec:typeA}
we apply our main theorem to intervals in the symmetric group $A_n$.
As a by-product a new Stirling number identity pops up
at the end of Section~\ref{sec:typeA}. In
Section~\ref{sec:typeB}
we apply our main theorem to the hyperoctahedral
group $B_n$. Finally, in Section~\ref{sec:open}
we discuss further research directions and
suggest some open problems.

\section{Rook polynomials}\label{sec:rook}
\begin{figure}
\begin{center}
\setlength{\unitlength}{0.5mm}
\resizebox{40mm}{!}{\includegraphics{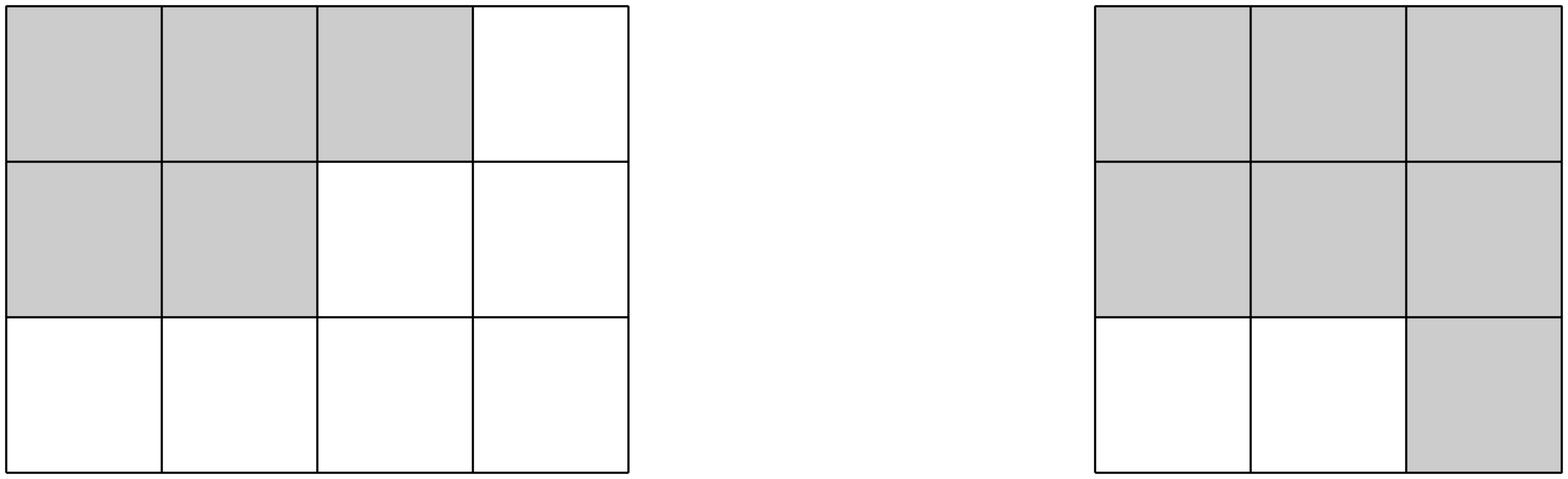}}
\caption[Square diagrams of Ferrers matrices.]
{Square diagrams of the left-aligned Ferrers matrix
$\lambda=\left(\begin{smallmatrix}
1 & 1 & 1 & 0 \\
1 & 1 & 0 & 0 \\
0 & 0 & 0 & 0
\end{smallmatrix}\right)$
and the right-aligned Ferrers matrix
$\mu=\left(\begin{smallmatrix}
1 & 1 & 1 \\
1 & 1 & 1 \\
0 & 0 & 1
\end{smallmatrix}\right)$. The row and column lengths
are given by the following table:
\[
\begin{array}{c|cccc}
i & r_i(\lambda) & c_i(\lambda) & r_i(\mu) & c_i(\mu) \\
\hline
1 & 3 & 2 & 3 & 2 \\
2 & 2 & 2 & 3 & 2 \\
3 & 0 & 1 & 1 & 3 \\
4 &   & 0 &   &
\end{array}
\]
}
\label{fig:ferrers}
\end{center}
\end{figure}

Let $A$ be a zero-one matrix and put rooks on some of the
one-entries of $A$. If no two rooks are in the same row or column
we have a {\em (non-taking) rook configuration on $A$}, and
we say that $A$ {\em covers} the rook configuration.
In the literature, $A$ is sometimes called a {\rm board}
and is often depicted by square diagrams like those in
Figure~\ref{fig:ferrers}. For convenience we will simultaneously
think of $A$ as the set of its one-entries, and write
for instance $(i,j)\in A$ if $A_{i,j}=1$ and use notation like
$A\cap B$.

Let $A^\updownarrow$ and $A^\crs$ denote reflecting the matrix upside
down respectively rotating it 180 degrees,
i.e.~$A^\updownarrow_{i,j}=A_{m-i+1,j}$ and
$A^\crs_{i,j}=A_{m-i+1,n-j+1}$ if $A$ is an $m\times n$ matrix.
Define $\pi^\updownarrow$ and $\pi^\crs$ similarly for rook
configurations $\pi$.

The number of rook configurations on $A$ with $k$ rooks is
called the {\em $k$th rook number of $A$} and is denoted
by $R^A_k$. Given a
nonnegative integer $n$, following Goldman
et al.~\cite{goldmanjoichiwhite1}
we define
{\em the $n$th rook polynomial of $A$} as
\[
\hat{R}^A_n(x)\DEF\sum_{k=0}^n R^A_{n-k} x(x-1)\dotsm(x-k+1).
\]
Note that $\hat{R}^A_n(0)=R^A_n$.

A zero-one matrix $\lambda$ is a {\em left-aligned (resp.~right-aligned)
Ferrers matrix} if
every one-entry has one-entries directly to the left (resp.~to the right)
and above it (provided these entries exist). The number of ones in the
$i$th row (resp.~column) of $\lambda$ is denoted by
$r_i(\lambda)$ (resp.~$c_i(\lambda)$).
Figure~\ref{fig:ferrers} shows an example.

From~\cite{goldmanjoichiwhite1} we have the following theorem.
\begin{theo}[Goldman et al.]\label{th:ferrers}
Let $\lambda$ be a right-aligned Ferrers matrix of size $m\times n$. Then
\[
\hat{R}^\lambda_n(x)=\prod_{j=1}^n(x+c_j(\lambda)-j+1).
\]
\end{theo}

Given a rook configuration $\mathcal{A}$ on A,
define the statistics $\inv_A(\mathcal{A})$
to be the number of (not necessarily positive) cells of $A$
with no rook weakly to the right in the same row or below
in the same column.
In the special case
where $A$ is an $n\times n$ matrix and $\mathcal{A}$ has
$n$ rooks, $\inv_A(\mathcal{A})$ becomes the number of
inversions of the permutation $\pi$ given by
$\pi(i)=j\ \Leftrightarrow\ (i,j)\in\mathcal{A}$,
where $i$ is the row index and $j$ is the column index.

Next, (almost) following Garsia and Remmel~\cite{garsiaremmel},
we define the {\em $k$th $q$-rook number of $A$} as
\[
R^A_k(q)=\sum_{\mathcal{A}} q^{\inv_A{\mathcal{A}}}
\]
where the sum is over all rook configurations on $A$ with
$k$ rooks. Given a nonnegative integer $n$,
the {\em $n$th $q$-rook polynomial of $A$} is defined as
\[
\hat{R}^A_n(x;q)
\DEF\sum_{k=0}^n R^A_{n-k}(q)[x]_q[x-1]_q\dotsm[x-k+1]_q.
\]
Here $[x]_q\DEF 1+q+q^2+\dotsb+q^{x-1}=(1-q^x)/(1-q)$
is the $q$-analogue of $x$.
Observe that putting $q=1$ yields the ordinary rook numbers and
polynomials.

Garsia and Remmel showed that
Theorem~\ref{th:ferrers} has a beautiful $q$-analogue:
\begin{theo}[Garsia, Remmel]
Let $A$ be a left-aligned Ferrers matrix of size $m\times n$. Then
\[
\hat{R}^A_n(x;q)=q^z\prod_{j=1}^n[x+c_j(A)+j-n]_q
\]
where $z$ is the number of zero-entries in $A$.
\end{theo}
Let $[n]!_q\DEF[1]_q[2]_q\dotsm[n]_q$.
\begin{coro}\label{co:square}
For the $n\times n$ square matrix $J^{n,n}$ with ones everywhere,
the $n$th $q$-rook number is
\[
R^{J^{n,n}}_n(q)=[n]!_q.
\]
\end{coro}

Let $T_n$ denote the $n\times n$ zero-one matrix with ones
on and above the secondary diagonal,
i.e.~$(T_n)_{i,j}=1\ \Leftrightarrow\ i\le n-j+1$.
In~\cite[p.~248]{garsiaremmel} it is proved that
\begin{equation}\label{eq:stirling}
R^{T_n}_k(q)=q^{\binom{n}{2}}S_{n+1,n+1-k}(q)
\end{equation}
where $S_{n,k}(q)$ is the {\em $q$-Stirling number} defined by
the recurrence
\[
S_{n+1,k}(q)=q^{k-1}S_{n,k-1}+[k]_qS_{n,k}(q),\ \ \ \mbox{for $0\le k\le n$}
\]
with the initial conditions $S_{0,0}(q)=1$ and $S_{n,k}=0$
for $k<0$ or $k>n$.

\section{Results}\label{sec:results}
A {\em skew Ferrers matrix} $\lambda/\mu$ is the difference
$\lambda-\mu$ between a Ferrers matrix $\lambda$ and an equally aligned
componentwise smaller Ferrers matrix $\mu$. If $\lambda$ and $\mu$
are left-aligned, then $\lambda/\mu$ is also said to be left-aligned,
and if $\lambda$ and $\mu$ are right-aligned, so is $\lambda/\mu$.

Let $\fs_n$ denote the symmetric group.
For any zero-one $n\times n$ matrix $A$, let $\fs(A)$ be the
set of rook configurations on $A$ with $n$ rooks.
We will identify such a rook configuration with
a permutation $\pi\in\fs_n$ so that $\pi(i)=j$
if and only if there is a rook at the square $(i,j)$,
where $i$ is the row index and $j$ is the column index.

\begin{figure}
\begin{center}
\setlength{\unitlength}{0.5mm}
\resizebox{50mm}{!}{\includegraphics{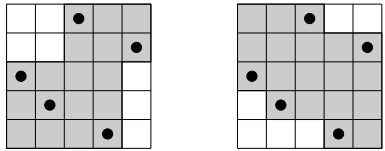}}
\caption{The shaded regions show the
left respectively right hull of the permutation 35124.}
\label{fig:hulls}
\end{center}
\end{figure}

For a permutation $\pi\in\fs_n$, let the {\em right (resp.~left) hull
$H_R(\pi)$ (resp.~$H_L(\pi)$) of $\pi$} be the smallest right-aligned
(resp.~left-aligned) skew Ferrers matrix that covers
$\pi$. Figure~\ref{fig:hulls} shows an example.

For the definition of Bruhat order and a general treatment
of Coxeter groups from a combinatorialist's viewpoint,
we refer to Bj\"orner and Brenti~\cite{bjornerbrenti}.

\begin{figure}
\begin{center}
\setlength{\unitlength}{0.5mm}
\resizebox{100mm}{!}{\includegraphics{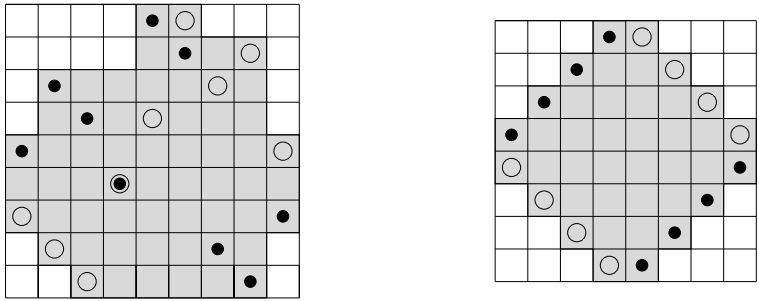}}
\caption[]{Left: The permutations $u=562314978$ (dots)
and $w=687594123$ (circles) in $\fs_9$
satisfy the pattern condition in Corollary~\ref{co:poincare}, so the
interval $[u,w]$ consists of precisely
the permutations that fit inside the shaded region $H_R(w)\cap H_L(u)$.
Right: $H_R(56781234)\cap H_L(43218765)$ is the Aztec diamond of order
4. (In fact it follows from part~(3) of Corollary~\ref{co:poincare} that
there are $2^n$ elements in the interval $[w^\updownarrow,w]$ in
$A_{2n-1}$, where $w=\max A_{2n-1}^{S\setminus\{s_n\}}$.)}
\label{fig:mainexample}
\end{center}
\end{figure}

Our main result is the following theorem and its corollary.
\begin{theo}\label{th:patterns}
$\fs(H_R(\pi))$ equals the lower Bruhat interval $[\id,\pi]$
in $\fs_n$
if and only if $\pi$ avoids the patterns 4231, 35142, 42513, and 351624.
\end{theo}
\begin{remark}
The permutations $\pi$ in the theorem are exactly the ones
such that the Schubert variety corresponding to the interval
$[\id,\pi]$ is {\em defined by inclusions} in the sense
of Gasharov and Reiner~\cite{gasharovreiner} according
to their Theorem~4.2.
\end{remark}
\begin{coro}\label{co:poincare}
  Let $u,w\in\fs_n$ and suppose
  $w$ and $u^\updownarrow$ both avoid the patterns
  4231, 35142, 42513, and 351624.
  Then the following holds.
  \begin{enumerate}
  \item
    $\fs(H_R(w)\cap H_L(u))$ equals the Bruhat
    interval $[u,w]$.
  \item
    The Poincar\'e polynomial $\poin_{[u,w]}(q)$ of $[u,w]$
    equals the $q$-rook number
    $R^{H_R(w)\cap H_L(u)}_n(q)$.
  \item
    In particular, the number of elements in $[u,w]$ equals
    the ordinary rook number $R^{H_R(w)\cap H_L(u)}_n$.
  \end{enumerate}
\end{coro}
\begin{proof}
  Once we observe that $H_L(u)=H_R(u^\updownarrow)^\updownarrow$
  and recall that flipping the
  rook configurations upside down
  is an antiautomorphism
  on the Bruhat order on $\fs_n$,
  the corollary follows directly from Theorem~\ref{th:patterns}.
\end{proof}
\begin{remark}
If $\pi$ is 312-avoiding then $H_R(\pi)$ is an ordinary Ferrers
matrix and Corollary~\ref{co:poincare} overlaps with
Theorem~33 in~\cite{ding} by Ding. In this case Ding
coined the name {\em partition variety} for
the Schubert variety corresponding to the
Bruhat interval $[\id,\pi]$ in $\fs_n$. Thus it would be
logical to coin the name {\em skew partition variety}
for a Schubert variety corresponding to an
interval $[\id,\pi]$ such that $[\id,\pi]=\fs(H_R(\pi))$.
\end{remark}
\noindent
Figure~\ref{fig:mainexample} shows two examples of the corollary.

Nontrivial application of the above result yields the
Poincar\'e polynomial of some particularly interesting
intervals in finite Weyl groups.

For a Coxeter system $(W,S)$ and a subset $J\subseteq S$ of
the generators, let $W_J$ denote the parabolic
subgroup generated by $J$. Each left coset $wW_J\in W/W_J$ has a unique
representative of minimal length, see~\cite[Cor.~2.4.5]{bjornerbrenti}.
The system of such minimal coset representatives is denoted
by $W^J$, and the Bruhat order on $W$ restricts to an order on $W^J$.

\begin{figure}
\begin{center}
\setlength{\unitlength}{0.5mm}
\begin{picture}(100,20)(0,-10)
\put(-30,-2){$A_n$}
\put(0,0){\line(1,0){50}}
\dashline[-2]{1}(45,0)(65,0)
\put(70,0){\line(1,0){30}}
\put(0,0){\circle*{4}}
\put(20,0){\circle*{4}}
\put(40,0){\circle*{4}}
\put(80,0){\circle*{4}}
\put(100,0){\circle*{4}}
\put(-4,-10){$s_1$}
\put(16,-10){$s_2$}
\put(36,-10){$s_3$}
\put(76,-10){$s_{n-1}$}
\put(96,-10){$s_n$}
\end{picture}\\
\begin{picture}(100,25)(0,-10)
\put(-30,-2){$B_n$}
\put(0,0){\line(1,0){50}}
\dashline[-2]{1}(45,0)(65,0)
\put(70,0){\line(1,0){30}}
\put(0,0){\circle*{4}}
\put(20,0){\circle*{4}}
\put(40,0){\circle*{4}}
\put(80,0){\circle*{4}}
\put(100,0){\circle*{4}}
\put(-4,-10){$s_0$}
\put(7,2){4}
\put(16,-10){$s_1$}
\put(36,-10){$s_2$}
\put(76,-10){$s_{n-2}$}
\put(96,-10){$s_{n-1}$}
\end{picture}
\caption{The Coxeter graphs of $A_n$ and $B_n$.}
\label{fig:coxetergraphs}
\end{center}
\end{figure}
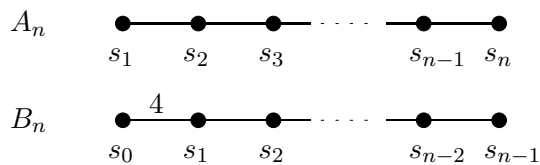

We will deal with two infinite families of finite Coxeter systems,
namely the symmetric groups
$A_n$ and the hyperoctahedral groups
$B_n$. Their Coxeter graphs are depicted in
Figure~\ref{fig:coxetergraphs}.

For type A we have the following result.
\begin{theo}\label{th:typeA}
Let $w$ be the maximal element of $A_{n-1}^{S\setminus\{s_k\}}$.
Then the Poincar\'e polynomial of the Bruhat interval $[\id,w]$
is
\[
\poin_{[\id,w]}(q)=q^{(n-k)k}\sum_{i=0}^k
S_{k+1,i+1}(1/q)\,S_{n-k+1,i+1}(1/q)\,[i]!_q^2\,q^i.
\]
\end{theo}
The special case where $n$ is even, $k=n/2$, and $q=1$
follows from Exercise~4.36 in Lov\'asz~\cite{lovasz}
once one knows that the set of permutations
he describes is a Bruhat interval (which
is Exercise~2.6 in~\cite{bjornerbrenti}). See also
Theorem~3 in Vesztergombi~\cite{vesztergombi}.

For type B the corresponding result looks like this:
\begin{theo}\label{th:typeB}
Let $w$ be the maximal element of $B_n^{S\setminus\{s_0\}}$.
Then the Poincar\'e polynomial of the Bruhat interval $[\id,w]$ is
\[
\poin_{[\id,w]}(q)=q^{\binom{n+1}{2}}
\sum_{i=0}^n S_{n+1,i+1}(1/q)\,[i]!_q.
\]
\end{theo}

We also present a recurrence relation for computing the
number of elements in the Bruhat interval $[\id,w]$ of $A_{n-1}$
where $w$ is any element in $A_{n-1}^{S\setminus\{s_k\}}$.
As a by-product we obtain the following
Stirling number identity which appears to be new.
\begin{theo}\label{th:bernoulli}
Let $w$ be the maximal element of $A_{n-1}^{S\setminus\{s_k\}}$.
Then the number of elements in the Bruhat interval $[\id,w]$
is
\[
\begin{split}
\poin_{[\id,w]}(1) & =\sum_{i=0}^k S_{k+1,i+1}\,S_{n-k+1,i+1}\,i!^2\\
& =(-1)^k\sum_{i=0}^k(-1)^i(i+1)^{n-k}i! S_{k,i} \\
& =B_n^{k-n}
\end{split}
\]
where $S_{n,k}$ are the Stirling numbers of the second kind,
and $B_n^k$ are the {\em poly-Bernoulli numbers}
defined by Kaneko~\cite{kaneko}.
\end{theo}
\begin{remark}
Kaneko's Theorem~2 says that $B_n^{-k}=B_k^{-n}$ for any $n,k\ge0$.
This is immediately evident
from our new formula for the poly-Bernoulli numbers.
\end{remark}
From Kaneko's work~\cite[p.~223]{kaneko} we can compute the
exponential bivariate generating function for $\poin_{[\id,w]}(1)$.
\[
\sum_{n=0}^\infty\sum_{k=0}^\infty
\poin_{[\id,w]}(1)\frac{x^n}{n!}\frac{y^{n-k}}{(n-k)!}
=\frac{e^{x+y}}{e^x+e^y-e^{x+y}}.
\]
The poly-Bernoulli numbers have the sequence number
A099594 in Sloane's On-Line Encyclopedia of Integer Sequences.

\section{Skew Ferrers matrices and Poincar\'e polynomials}
\label{sec:connection}
In this section we make a connection between
Poincar\'e polynomials and rook polynomials, and
prove Theorem~\ref{th:patterns}.
\begin{prop}
If $\lambda/\mu$ is a right-aligned skew Ferrers matrix of size $n\times n$,
then $\fs(\lambda/\mu)$ is an order ideal in the Bruhat order of $\fs_n$.
\end{prop}
\begin{proof}
The Bruhat order is the transitive closure of the Bruhat graph
whose edges correspond to transpositions
(see~\cite[Sec.~2.1]{bjornerbrenti}). Thus it suffices to
show that we cannot leave $\lambda/\mu$ by a transposition
going down in the Bruhat order. In other words,
if $\pi$ is a rook configuration on $\lambda/\mu$ with $n$ rooks,
and $\pi_i>\pi_{i'}$ with $i<i'$, then exchanging rows
$i$ and $i'$ yields a rook configuration which is covered
by $\lambda/\mu$. This is obviously true, as we can see in
Figure~\ref{fig:orderideal}.
\begin{figure}
\begin{center}
\setlength{\unitlength}{0.5mm}
\resizebox{50mm}{!}{\includegraphics{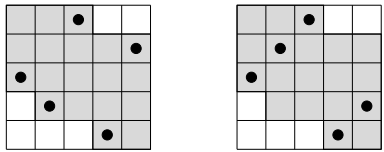}}
\caption{The permutation 35124 to the left becomes 32154 to the right
after exchanging rows 2 and 4. We do not leave the shaded region
$\lambda/\mu$ by this operation.}
\label{fig:orderideal}
\end{center}
\end{figure}
\end{proof}

For $\pi\in\fs_n$ and $i,j\in[n]=\{1,2,\dotsc,n\}$, let
\[
\pi[i,j]\DEF\{a\in[i]\ :\ \pi(a)\ge j\}.
\]
In other words
$\pi[i,j]$ is the number of rooks weakly north-east
of the square $(i,j)$.
The following criterion for comparing two permutations
with respect to the Bruhat order is well-known
(see e.g.~\cite[Th.~2.1.5]{bjornerbrenti}).
\begin{lemma}\label{lm:bruhat}
Let $\pi,\rho\in\fs_n$. Then $\pi\le\rho$ if and only if
$\pi[i,j]\le\rho[i,j]$ for all $i,j\in[n]$.
\end{lemma}

Theorem~\ref{th:patterns} completely characterises the
interesting cases where $\fs(\lambda/\mu)$
is a lower Bruhat interval $[\id,\pi]$. Now we are
ready for the proof.
\begin{proof}[\bf Proof of Theorem~\ref{th:patterns}]
We begin with the ``only if'' direction which is the easier one.
For each of the four forbidden patterns we will do the following:
First we suppose $\pi$ contains the pattern. Then
we move some of the rooks that constitute the pattern to
new positions, and call the resulting rook configuration $\rho$.
This $\rho$ is seen to be covered by
$H_R(\pi)$ while $\rho\not\le\pi$ in Bruhat order,
and we conclude that $\pi$ is not uniquely maximal in $H_R(\pi)$.

Suppose $\pi$ contains the pattern 4231 so that there
are rooks $(i_1,j_4)$, $(i_2,j_2)$, $(i_3,j_3)$, and
$(i_4,j_1)$ with $i_1<i_2<i_3<i_4$ and $j_1<j_2<j_3<j_4$.
Move the rooks $(i_2,j_2)$ and $(i_3,j_3)$ to the positions
$(i_2,j_3)$ and $(i_3,j_2)$ and call the resulting
rook configuration $\rho$. Then $\rho$
is covered by $H_R(\pi)$ and
$\rho>\pi$ in Bruhat order so $\pi$ is not maximal in $H_R(\pi)$.

Note that the rooks outside the pattern turned out to be irrelevant for
the discussion. In fact we could have supposed $\pi$
was equal to the pattern 4231 and then simply defined $\rho=4321$.
This observation applies to the remaining three patterns as well,
and thus the ``only if'' part of the proof can be
written as a table that associates a $\rho$ to each pattern $\pi$:
\begin{center}
\begin{tabular}{c|c}
$\pi$ & $\rho$ \\
\hline 4231 & 4321 \\
35142 & 15432 \\
42513 & 43215 \\
351624 & 154326
\end{tabular}
\end{center}
Figure~\ref{fig:patterns} illustrates the table and makes it evident
that $\rho$ is covered by $H_R(\pi)$ in each case. That
$\rho\not\le\pi$ in Bruhat order can be checked easily
using Lemma~\ref{lm:bruhat}.
\begin{figure}
\begin{center}
\setlength{\unitlength}{0.5mm}
\resizebox{100mm}{!}{\includegraphics{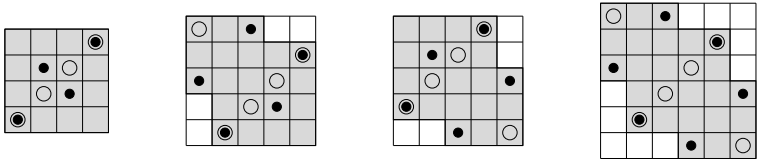}}
\caption{The dots show the rook configuration $\pi$ in the four cases
4231, 35142, 42513, and 351624.
The shaded squares show the right hull $H_R(\pi)$,
and the circles show $\rho$.}
\label{fig:patterns}
\end{center}
\end{figure}




Now it is time to prove the difficult ``if'' direction.
Suppose $\fs(H_R(\pi))\ne[\id,\pi]$ so that there is a
$\rho\in\fs(H_R(\pi))$ with $\rho\not\le\pi$. Our goal is to show
that $\pi$ contains some of the four forbidden patterns.

Let the rooks of $\pi$ and $\rho$ be black and white, respectively.
(Observe that some squares may contain both a black and a white rook.)
Order the squares $[n]^2$ partially
so that $(i,j)\le(i',j')$ if $i\le i'$ and $j\ge j'$,
i.e.~the north-east corner $(1,n)$ is the minimal square
of $[n]^2$.

Let $L$ be the set of squares $(i,j)$ with $\rho[i,j]>\pi[i,j]$ and
no black rook weakly to the right of $(i,j)$ in row $i$ or
above $(i,j)$ in column $j$.
First we show that $L$ is not empty.

Since $\rho\not\le\pi$, by Lemma~\ref{lm:bruhat}
there is a square $(i,j)\in[n]^2$
such that $\rho[i,j]>\pi[i,j]$. Let $(i_{\rm min},j_{\rm min})$ be a minimal
square with this property. Then there is no black rook weakly
to the right of $(i_{\rm min},j_{\rm min})$ in row $i_{\rm min}$,
for if that were the case the smaller square
$(i_{\rm min}-1,j_{\rm min})$ would have the property
$\rho[i_{\rm min}-1,j_{\rm min}]>\pi[i_{\rm min}-1,j_{\rm min}]$ as well.
Analogously,
there is no black rook weakly above $(i_{\rm min},j_{\rm min})$
in column $j_{\rm min}$. Thus $(i_{\rm min},j_{\rm min})$ belongs to $L$.

Now we can let $(i_{\rm max},j_{\rm max})$
be a maximal square in $L$. Since
$\rho[i_{\rm max},j_{\rm max}]>\pi[i_{\rm max},j_{\rm max}]$
we have $i_{\rm max}<n$ and $j_{\rm max}>1$.
There must be a black rook weakly to the right
of $(i_{\rm max}+1,j_{\rm max})$ in row $i_{\rm max}+1$ because otherwise
the greater square $(i_{\rm max}+1,j_{\rm max})$ would belong to $L$.
By an analogous argument, there is a black rook weakly above
$(i_{\rm max},j_{\rm max}-1)$ in column $j_{\rm max}-1$.
We have the situation depicted in Figure~\ref{fig:situation}.
\begin{figure}
\begin{center}
  \setlength{\unitlength}{0.3mm}
    \begin{picture}(130,130)(-10,0)
      \path(10,0)(110,0)(110,100)(10,100)(10,0)
      \texture{cccccccc 0 0 0 0 0 0 0
        cccccccc 0 0 0 0 0 0 0
        cccccccc 0 0 0 0 0 0 0
        cccccccc 0 0 0 0 0 0 0}
      \shade\path(10,50)(50,50)(50,0)(40,0)(40,40)(10,40)(10,50)
      \shade\path(50,50)(110,50)(110,60)(60,60)(60,100)(50,100)(50,50)
            
      \put(-28,53){$\scriptstyle i_{\rm max}$}
      \put(-5,55){\vector(1,0){12}}
      \put(48,120){$\scriptstyle j_{\rm max}$}
      \put(55,115){\vector(0,-1){12}}
    \end{picture}
\end{center}
\caption{The shaded region contains no black rook.}
\label{fig:situation}
\end{figure}
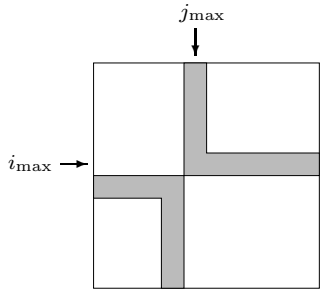

Since there are more white than black rooks inside the rectangle
$R=[1,i_{\rm max}]\times[j_{\rm max},n]$ there must also
be more white than black rooks inside the diagonally opposite
rectangle $R'=[i_{\rm max}+1,n]\times[1,j_{\rm max}-1]$.
In particular there is at least one white rook inside $R$
and at least one white rook inside $R'$. Since $\rho$ is
covered by $H_R(\pi)$ it follows that there is
a black rook $(i,j)$ inside $R$ and a black rook $(i',j')$
inside $R'$; choose $(i',j')$ minimal in $R'$.
Call $(i,j)$ and $(i',j')$ {\em the witnesses}.
Now the situation is exactly as in the proof of
Theorem 4.2 in~\cite{gasharovreiner} by Gasharov and Reiner.
The remaining part of the proof will
essentially be a copy of their arguments.

We show that at least one of the four forbidden patterns will appear,
depending on whether the rectangle $[i,i']\times[j',j]$
contains a black rook strictly to the left of column $j_{\rm max}$, and
a black rook strictly below row $i_{\rm max}$.
If one can find
\begin{enumerate}
\item
  both, then combining these with the two witnesses produces
  the pattern 4231 in $\pi$.
  (Look at Figure~\ref{fig:1234} for illustrations.)
\item
  the former but not the latter, then combining the two witnesses with
  the former and with the black rooks in column $j_{\rm max}$
  and in row $i_{\rm max}+1$ produces the pattern 42513.
\item
  the latter but not the former, then combining the two witnesses with the
  latter and with the black rooks in column $j_{\rm max}-1$
  and in row $i_{\rm max}$ produces the pattern 35142.
\item
  neither, then combining the two witnesses with the black rooks
  in column $j_{\rm max}-1$ and $j_{\rm max}$ and in row
  $i_{\rm max}$ and $i_{\rm max}+1$ produces the pattern 351624.
\end{enumerate}

  \newsavebox{\situation}
  \savebox{\situation}(0,0)[lb]{\setlength{\unitlength}{0.3mm}%
    \begin{picture}(160,140)(-10,0)
      \path(10,0)(110,0)(110,100)(10,100)(10,0)
      \texture{cccccccc 0 0 0 0 0 0 0
        cccccccc 0 0 0 0 0 0 0
        cccccccc 0 0 0 0 0 0 0
        cccccccc 0 0 0 0 0 0 0}
      \shade\path(10,50)(50,50)(50,0)(40,0)(40,40)(10,40)(10,50)
      \shade\path(25,25)(40,25)(40,40)(25,40)(25,25)
      \shade\path(50,50)(110,50)(110,60)(60,60)(60,100)(50,100)(50,50)
      
      \put(25,25){\circle*{6}}
      \put(85,75){\circle*{6}}
      
      \path(25,25)(85,25)(85,75)(25,75)(25,25)
      
      \put(-28,53){$\scriptstyle i_{\rm max}$}
      \put(-5,55){\vector(1,0){12}}
      \put(48,120){$\scriptstyle j_{\rm max}$}
      \put(55,115){\vector(0,-1){12}}
      
      \put(-15,73){$\scriptstyle i$}
      \put(-5,75){\vector(1,0){12}}
      \put(83,120){$\scriptstyle j$}
      \put(85,115){\vector(0,-1){12}}
      
      \put(-15,23){$\scriptstyle i'$}
      \put(-5,25){\vector(1,0){12}}
      \put(23,120){$\scriptstyle j'$}
      \put(25,115){\vector(0,-1){12}}
      
    \end{picture}}
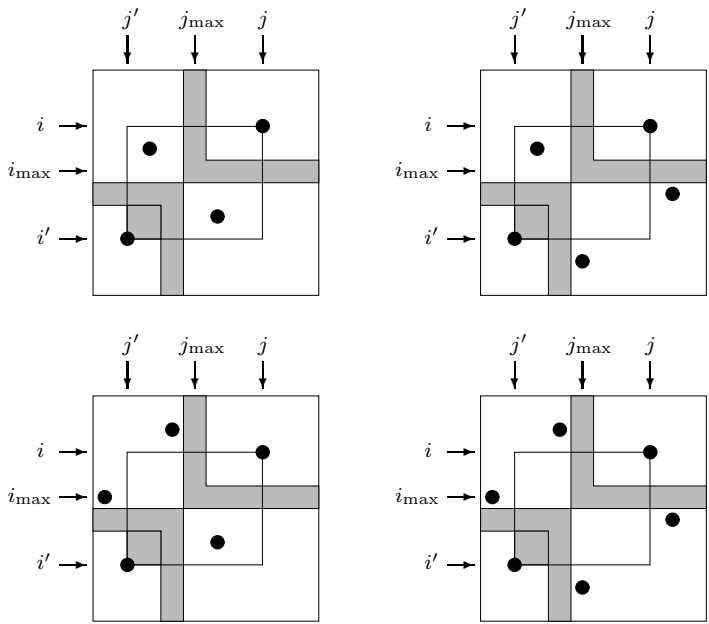
\begin{figure}
\begin{center}
  \setlength{\unitlength}{0.3mm}
  \begin{tabular}{cc}
  \begin{picture}(160,140)(-10,0)
  \put(0,0){\usebox{\situation}}
      \put(45,65){\circle*{6}}
      \put(75,35){\circle*{6}}

\end{picture} &
  \begin{picture}(160,140)(-10,0)
  \put(0,0){\usebox{\situation}}
      \put(45,65){\circle*{6}}
      \put(65,15){\circle*{6}}
      \put(105,45){\circle*{6}}

\end{picture} \\
  \begin{picture}(160,140)(-10,0)
  \put(0,0){\usebox{\situation}}
      \put(25,55){\circle*{6}}
      \put(55,85){\circle*{6}}
      \put(75,35){\circle*{6}}

\end{picture} &
  \begin{picture}(160,140)(-10,0)
  \put(0,0){\usebox{\situation}}
      \put(25,55){\circle*{6}}
      \put(55,85){\circle*{6}}
      \put(65,15){\circle*{6}}
      \put(105,45){\circle*{6}}

\end{picture}
\end{tabular}
\caption{The four cases of the ``if'' part of the proof
of Theorem~\ref{th:patterns}. The shaded regions contain no black rooks.}
\label{fig:1234}
\end{center}
\end{figure}

\end{proof}

\section{Poincar\'e polynomials of $A_n$}\label{sec:typeAq}
In this section we apply Theorem~\ref{th:patterns} to
the lower Bruhat interval $[\id,w]$ of the symmetric group
$A_n$ where $w$ is the maximal minimal coset representative
$w=\max A_n^{S\setminus\{s_k\}}$. In the end we obtain the
simple formula of Theorem~\ref{th:typeA}.

Let $J^{m,n}$ denote the $m\times n$ matrix with all entries equal
to one.
The following is a $q$-analogue of the corollary to Theorem~1
in~\cite{cheonhwangsong}.
\begin{prop}\label{pr:sharp}
Let $A$ and $B$ be zero-one matrices of sizes
$m\times m$ and $n\times n$, respectively.
The block matrix
\[
B\kors A\DEF\begin{pmatrix}
B & J^{n,m} \\
J^{m,n} & A
\end{pmatrix}
\]
has the $(m+n)$th $q$-rook number
\[
R^{B\kors A}_{m+n}(q)
=\sum_{i=0}^{\min(m,n)} R^A_{m-i}(q)\,
R^{B^\crs}_{n-i}(q)\,[i]!_q^2\,q^{-i^2}.
\]
\end{prop}
\begin{proof}
It is easy to see that
each configuration $\pi$ of $m+n$ rooks on $B\kors A$ is
chosen uniquely by the following procedure:
\begin{itemize}
\item
First, choose a nonnegative integer $i$.
\item
Then choose a configuration $\mathcal{A}$ of
$m-i$ rooks on $A$ and a configuration $\mathcal{B}$ of
$n-i$ rooks on $B$. Together $\mathcal{A}$ and $\mathcal{B}$
form a configuration of $m+n-2i$ rooks on $\mat{B}{0}{0}{A}$.
\item
Let $X$ be the $i\times i$ submatrix consisting of
the remaining free one-entries of $\mat{0}{J^{n,m}}{0}{0}$, i.e.~the
one-entries whose row and column have no rook in $\mathcal{A}$
or $\mathcal{B}$. Similarly, let $Y$ be the
$i\times i$ submatrix consisting of
the remaining free one-entries of $\mat{0}{0}{J^{m,n}}{0}$.
Now choose a configuration $\mathcal{X}$ of
$i$ rooks on $X$ and a configuration $\mathcal{Y}$
of $i$ rooks on $Y$.
\end{itemize}
Let $\Inv(\pi)$ be the set of inversions of $\pi$, i.e.~pairs $(r,r')$ of rooks
such that $r$ is strictly north-east of $r'$.
The number $\inv_{A}(\mathcal{A})$ counts the
cells in $A$ which have no rooks to the right or below.
This equals the number of inversions $(r,r')$
such that $r$ belongs to
$A$ or $J^{n,m}$ and $r'$ belongs to $A$ or $J^{m,n}$:
\[\inv_{A}(\mathcal{A})=\abs{\{(r,r')\in\Inv(\pi)
\ :\ r\in\mat{0}{J^{n,m}}{0}{A},\ r'\in\mat{0}{0}{J^{m,n}}{A}\}}.
\]
Similarly, $\inv_{B^\crs}(\mathcal{B}^\crs)$ counts the cells
in $B$ which have no rooks to the left or above, so
\[
\inv_{B^\crs}(\mathcal{B}^\crs)=\abs{\{(r,r')\in\Inv(\pi)
\ :\ r\in\mat{B}{J^{n,m}}{0}{0},\ r'\in\mat{B}{0}{J^{m,n}}{0}\}}.
\]
We also have
\[
\inv_X(\mathcal{X})=\abs{\{(r,r')\in\Inv(\pi)
\ :\ r,r'\in\mat{0}{0}{J^{m,n}}{0}\}}.
\]
and
\[
\inv_Y(\mathcal{Y})=\abs{\{(r,r')\in\Inv(\pi)
\ :\ r,r'\in\mat{0}{J^{n,m}}{0}{0}\}}.
\]
Putting the above equations together yields
\begin{eqnarray}\label{eq:inv}
&&\inv_{A}(\mathcal{A})+\inv_{B^\crs}(\mathcal{B^\crs})+\inv_X(\mathcal{X})
+\inv_Y(\mathcal{Y})\nonumber\\
&=&\inv(\pi)+\abs{\{(r,r')\in\Inv(\pi)\ :\ r\in \mat{0}{J^{n,m}}{0}{0},
\ r'\in\mat{0}{0}{J^{m,n}}{0}\}}\\
&=&\inv(\pi)+i^2\nonumber
\end{eqnarray}
where $\inv(\pi)=\abs{\Inv(\pi)}$.
Now we exponentiate and
sum over all permutations $\pi$ which can be constructed
by the procedure above:
\[
\sum_\pi q^{\inv(\pi)}=\sum_{i=0}^{\min(m,n)} q^{-i^2}
R^{A}_{m-i}(q)R^{B^\crs}_{n-i}(q)R^X_i(q)R^Y_i(q).
\]
By Corollary~\ref{co:square}, $R^X_i(q)=R^Y_i(q)=[i]!_q$.
\end{proof}

\begin{proof}[Proof of Theorem~\ref{th:typeA}]
A Coxeter system of type $(A_{n-1},S=\{s_1,s_2,\dotsc,s_{n-1}\})$
(see Figure~\ref{fig:coxetergraphs})
is isomorphic to the symmetric group
$\fs_n$ with the adjacent transpositions $s_i=(i\leftrightarrow i+1)$
as generators.
A permutation $w\in\fs_n$ can be represented by a rook configuration
on $J^{n,n}$ with $n$ rooks, so that $w(i)=j$ if and only if
there is a rook in the cell $(i,j)$.

Let $w$ be the maximal element in
$A_{n-1}^{S\setminus\{s_k\}}$,
i.e.~$\bigl(w(1),w(2),\dotsc,w(n)\bigr)
=\bigl(n-k+1,n-k+2,\dotsc,n,1,2,\dotsc,n-k\bigr)$.
Then we have
$H_R(w)=(T_{n-k}^\crs\kors T_k)^\updownarrow$
and hence
\[
\poin_{[\id,w]}(q)=R^{(T_{n-k}^\crs\kors T_k)^\updownarrow}_n(q)
=q^{\binom{n}{2}}R^{T_{n-k}^\crs\kors T_k}_n(1/q)
\]
which by Proposition~\ref{pr:sharp} equals
\[
q^{\binom{n}{2}}\sum_{i=0}^{\min(k,n-k)}
R^{T_{n-k}}_{n-k-i}(1/q)
R^{T_k}_{k-i}(1/q)
[i]!_{1/q}^2q^{i^2}.
\]
Using Equation~\eqref{eq:stirling} we obtain
\[
P(q)=
q^{\binom{n}{2}}\sum_{i=0}^{\min(k,n-k)}
q^{-{\binom{n-k}{2}}}S_{n-k+1,i+1}(1/q)
q^{-{\binom{k}{2}}}S_{k+1,i+1}(1/q)
[i]!_{1/q}^2q^{i^2}.
\]
Since $[i]!_{1/q}=[i]!_q\cdot q^{-{\binom{i}{2}}}$ we are done.
\end{proof}

\section{The number of elements in a lower interval of $A_n$}\label{sec:typeA}
For a general minimal coset representative
$w\in A_{n-1}^{S\setminus\{s_k\}}$, it seems
very hard to compute the complete Poincar\'e polynomial.
In this section we will solve the easier problem
to determine $\poin_{[\id,w]}(1)$, i.e.~the
number of elements of $[\id,w]$.
We obtain a recurrence relation that allows us
to count the elements in polynomial
time. In the special case when $w$ is the maximal element
in $A_{n-1}^{S\setminus\{s_k\}}$ this method results in a
formula different from what we get if we put $q=1$ in
Theorem~\ref{th:typeA}. From this, rather unexpectedly, we obtain
an identity of Stirling numbers that we have not seen in the literature.

\begin{figure}
\begin{center}
\setlength{\unitlength}{0.5mm}
\resizebox{100mm}{!}{\includegraphics{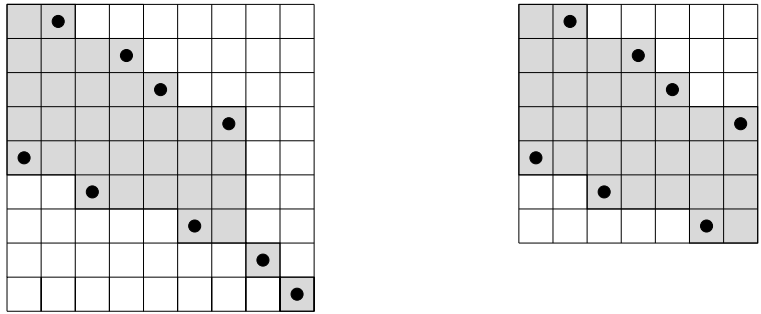}}
\caption{If $w(k)<n$ as in the left example, we may instead study the
smaller example to the right. They have isomorphic lower intervals.}
\label{fig:general}
\end{center}
\end{figure}

The set of minimal coset representatives
$A_{n-1}^{S\setminus\{s_k\}}$ consists of the permutations
$w\in\fs_n$
with $w(1)<w(2)<\dotsb<w(k)$ and $w(k+1)<w(k+2)<\dotsb<w(n)$.
Such a permutation clearly avoids the patterns in Corollary~\ref{co:poincare}
so the number of elements in the Bruhat interval $[\id,w]$
is given by the $n$th rook number $R^{H_R(w)}_n$.
Fortunately $H_R(w)$ has a simple structure.
If $w(k)<n$ then, as can be seen in Figure~\ref{fig:general},
\[
\bigl(w(w(k)+1),w(w(k)+2),\dotsc,w(n)\bigr)
=\bigl(n-w(k)+1,n-w(k)+2,\dotsc,n\bigr),
\]
so the interval $[\id,w]$ is isomorphic (as a poset) to the interval
$[\id,w']$ in $A_{w(k)-1}$, where $w'(i)=w(i)$ for $i=1,2,\dotsc,w(k)$.
Thus we may
assume that $w(k)=n$. Then
$H_R(w)=\lambda/\mu$, where $\lambda$ and $\mu$ are right-aligned
Ferrers matrices with row lengths
\begin{align*}
r_i(\lambda) & =\begin{cases}
n & \text{if $1\le i\le k$}, \\
n-w(i)+1 & \text{if $k+1\le i\le n$},
\end{cases} \\
r_i(\mu) & =\begin{cases}
n-w(i) & \text{if $1\le i\le k$}, \\
0 & \text{if $k+1\le i\le n$}.
\end{cases}
\end{align*}

For $1\le i\le k$, let $P_i$ be
the $n\times n$ zero-one matrix with ones in the cells
$(i,1),(i,2),\dotsc,(i,w(i))$.
It is easy
to see that a rook configuration with $n$ rooks is covered
by $\lambda/\mu$ if and only if it is covered by $\lambda$
and not by any $\lambda-P_i$. Thus, by the principle of inclusion-exclusion
we get
\[
R^{\lambda/\mu}_n=\sum_{I\subseteq[k]}(-1)^{\abs{I}}
R^{\lambda-\cup_{i\in I}P_i}_n.
\]
By a suitable permutation of the rows, the matrix
$\lambda-\cup_{i\in I}P_i$ can be transformed to a Ferrers matrix $\nu$
with column lengths
\[
c_j(\nu) = c_j(\lambda)-\abs{\{i\in I\ :\ w(i)\ge j\}}
=c_j(\lambda)-\abs{w(I)\cap[j,n]},
\]
where $w(I)=\{w(i)\ :\ i\in I\}$
denotes the image of $I$ under $w$.
Theorem~\ref{th:ferrers} with $x=0$ gives
\[
R^{\lambda/\mu}_n=\hat{R}^{\lambda/\mu}_n(0)
=\sum_{J\subseteq w([k])}(-1)^{\abs{J}}
\prod_{j=1}^n(c_j(\lambda)-\abs{J\cap[j,n]}-j+1).
\]
As we will see in a moment, this expression can be computed
efficiently by dynamic programming.

For $1\le a\le n$ and $0\le b\le n$, let
\[
f(a,b)=\sum_{J\in{\binom{w([k])\cap[a,n]}{b}}}\prod_{j=a}^n
(c_j(\lambda)-\abs{J\cap[j,n]}-j+1)
\]
where $\binom{w([k])\cap[a,n]}{b}$ denotes the set of subsets of
$w([k])\cap[a,n]$ of size $b$. Also put $f(a,b)=0$ if $b<0$.
It is straightforward to verify the following recurrence relation:
\begin{equation}\label{eq:frec}
\begin{cases}
f(a,b)=\bigl(c_a(\lambda)-a-b+1\bigr)f(a+1,b)
& \text{if $n>a\not\in w([k])$},\\
f(a,b)=\bigl(c_a(\lambda)-a-b+1\bigr)\bigl(f(a+1,b)+f(a+1,b-1)\bigr)
& \text{if $n>a\in w([k])$},\\
f(n,b)=\delta_{b,0}.&
\end{cases}
\end{equation}
Here $\delta_{b,0}$ is Dirac's $\delta$-function which is 1 if $b=0$
and 0 otherwise. Since
\begin{equation}\label{eq:fsumma}
R^{\lambda/\mu}_n=\sum_{b=0}^k (-1)^b f(1,b)
\end{equation}
the number of elements in $[\id,w]$ is computable
in polynomial time.

A special application of the method above admits us to
prove our by-product Theorem~\ref{th:bernoulli}.
\begin{proof}[Proof of Theorem~\ref{th:bernoulli}]
Consider the case when $w$ is the maximal element in
$A_{n-1}^{S\setminus\{s_k\}}$,
i.e.~$\bigl(w(1),w(2),\dotsc,w(n)\bigr)
=\bigl(n-k+1,n-k+2,\dotsc,n,1,2,\dotsc,n-k\bigr)$.
Then $c_a(\lambda)=k+a$ if $a\le n-k$ and $c_a(\lambda)=n$
if $a\ge n-k+1$, so the recurrence~\eqref{eq:frec} becomes
\[
\begin{cases}
f(a,b)=(k-b+1)f(a+1,b) & \text{if $a\le n-k$}, \\
f(a,b)=(n-a-b+1)\bigl(f(a+1,b)+f(a+1,b-1)\bigr)
& \text{if $n-k+1\le a\le n-1$}. \\
f(n,b)=\delta_{b,0} &
\end{cases}
\]
Iteration of the first line of the recurrence yields
\begin{equation}\label{eq:smallb}
f(1,b)=(k-b+1)^{n-k}f(n-k+1,b).
\end{equation}
Putting $g(\alpha,\beta)=f(n-\alpha+1,\alpha-\beta)$ for
$1\le\alpha\le k$, our recurrence transforms to
\[
\begin{cases}
  g(\alpha,\beta)=\beta\cdot\bigl(g(\alpha-1,\beta)+g(\alpha-1,\beta-1)\bigr)
  & \text{if $2\le \alpha\le k$}. \\
  g(1,\beta)=\delta_{\beta,1} &
\end{cases}
\]
We recognise this as the recurrence for $\beta!S_{\alpha,\beta}$
where $S_{\alpha,\beta}$ are Stirling numbers
of the second kind; thus
$f(a,b)=(n-a-b+1)!S_{n-a+1,n-a-b+1}$ for $n-k+1\le a\le n$.
Combining this with Equation~\eqref{eq:smallb} and plugging the
result into Equation~\eqref{eq:fsumma}, we obtain
\[
R^{\lambda/\mu}_n=\sum_{b=0}^k (-1)^b(k-b+1)^{n-k}(k-b)!S_{k,k-b}
\]
which also can be written as
\[
(-1)^k\sum_{i=0}^k(-1)^i(i+1)^{n-k}i! S_{k,i}.
\]
This happens to be the formula for the poly-Bernoulli number
$B_n^{k-n}$ defined by Kaneko~\cite{kaneko}.
\end{proof}

\section{Type B}\label{sec:typeB}
In this section we compute the Poincar\'e polynomial
of the lower Bruhat interval $[\id,w]$ in the hyperoctahedral group
$B_n$, where
$w$ is the maximal minimal coset representative,
$w=\max B_n^{S\setminus\{s_0\}}$.

We will represent $B_n$ combinatorially
by the set $\fs_n^B\DEF\{\pi\in\fs_{2n}\ :\ \pi^\crs=\pi\}$
of rotationally symmetric maximal
rook configurations on $J^{2n,2n}$, see~\cite[Ch.~8]{bjornerbrenti}.
In this representation $\pi\le\rho$ in Bruhat order on $B_n$
if and only if $\pi\le\rho$ as elements of $\fs_{2n}$
(\cite[Cor.~8.1.9]{bjornerbrenti}). The rank of $\pi$ is
\begin{equation}\label{eq:rankB}
\ell(\pi)=\bigl(\inv(\pi)+\,\nega(\pi)\bigr)/2
\end{equation}
where $\inv(\pi)$
is the usual inversion number of $\pi$ as an element of $\fs_{2n}$,
and $\nega(\pi)\DEF\abs{\{i\in[n+1,2n]\ :\ \pi(i)\le n\}}$,
see~\cite[Ch.~8, Exercise~2]{bjornerbrenti}.

For a zero-one matrix $A$ of size $2n\times2n$, let
\[
\RB^A(q,t)\DEF\sum_{\pi\in\fs^B_n\cap\fs(A)} q^{\inv(\pi)}t^{\nega(\pi)}.
\]
\begin{prop}\label{pr:sharpB}
Let $A$ be a zero-one matrix of size $n\times n$.
Then
\[
\RB^{A^\crs\kors A}(q,t)
=\sum_{i=0}^n R^{A}_{n-i}(q^2)\,
\,[i]!_{q^2}\,q^{-i^2}t^i.
\]
\end{prop}
\begin{proof}
The proof is almost identical to the proof of
Proposition~\ref{pr:sharp}.

When $m=n$, it is easy to see
that the permutation $\pi$ constructed by the
procedure in the proof of Proposition~\ref{pr:sharp}
is rotationally symmetric if and only if
${\mathcal B}={\mathcal A}^\crs$ and
${\mathcal Y}={\mathcal X}^\crs$.
Thus, putting $m:=n$, $B:=A^\crs$, $Y:=X^\crs=X$,
${\mathcal B}:={\mathcal A}^\crs$, and
${\mathcal Y}:={\mathcal X}^\crs$ into
Equation~\eqref{eq:inv} and using the identity
$(A^\crs)^\crs=A$, we obtain
\[
2\left(\inv_A(\mathcal{A})+\inv_X(\mathcal{X})\right)=\inv(\pi)+i^2.
\]
Obviously, $\nega(\pi)=i$.
Exponentiation and
summation over all rotationally symmetric permutations $\pi$
on $A^\crs\kors A$ yields
\[
\sum_\pi q^{\inv(\pi)}t^{\nega(\pi)}=\sum_{i=0}^n q^{-i^2}
R^{A}_{n-i}(q^2)R^X_i(q^2)t^i.
\]
By Corollary~\ref{co:square} $R^X_i(q^2)=[i]!_{q^2}$.
\end{proof}

Now we are ready for the proof of Theorem~\ref{th:typeB}.
\begin{proof}[Proof of Theorem~\ref{th:typeB}]
From Theorem~\ref{th:patterns} and
Equation~\eqref{eq:rankB} we obtain
\[
\poin_{[\id,w]}(q)=\sum_{u\in[\id,w]}q^{\ell(u)}
=\RB^{H_R(w)}(q^{1/2},q^{1/2}).
\]
It is easy to see that $H_R(w)=(T_n^\crs\kors T_n)^\updownarrow$.
Hence
\[
\poin_{[\id,w]}(q)=\RB^{(T_n^\crs\kors T_n)^\updownarrow}(q^{1/2},q^{1/2})
=q^{n^2}\RB^{(T_n^\crs\kors T_n)}(q^{-1/2},q^{-1/2})
\]
where we have used the fact that
\[
\RB^{A^\updownarrow}(q,t)=q^{\binom{2n}{2}}t^n\RB^{A}(q^{-1},t^{-1})
\]
for any $2n\times2n$ zero-one matrix $A$.
By Proposition~\ref{pr:sharpB} we can now compute
\[
\poin_{[\id,w]}(q)=
q^{n^2}\sum_{i=0}^n R^{T_n}_{n-i}(1/q)[i]!_{1/q}q^{\binom{i}{2}}.
\]
Using Equation~\eqref{eq:stirling} and the identity
$[i]!_{1/q}q^{\binom{i}{2}}=[i]!_q$, we finally obtain
\[
\poin_{[\id,w]}(q)=
q^{\binom{n+1}{2}}\sum_{i=0}^n S_{n+1,i+1}(1/q)[i]!_q.
\]
\end{proof}

\section{Open problems}\label{sec:open}
Perhaps the reason we still do not know much about the Poincar\'e
polynomials of Bruhat intervals after several decades of research
in the area is the lack of natural methods to attack it.
We hope the framework and the tools presented here will
make the problem more accessible, and we would like to suggest
a number of interesting open questions.
\begin{itemize}
\item
What is the Poincar\'e polynomial
$\poin_{[\id,w]}(q)$ in the even-signed permutation
group $D_n$ if 
$w=\max D_n^{S\setminus\{s_0\}}$ is the maximal minimal
representative in the quotient modulo a maximal parabolic subgroup
isomorphic to $A_{n-1}$?
\item
What is $\poin_{[\id,w]}(q)$ in the affine group $\tilde{A}_n$?
\item
Are there formulas for the generalized Eulerian polynomial
$\sum_{v\in[\id,w]}t^{d(v)}$ or even for the
bivariate generating function
$\sum_{v\in[\id,w]}t^{d(v)}q^{\ell(v)}$,
where $d(v)=\abs{\{s\in S\ :\ \ell(vs)<\ell(v)\}}$
is the {\em descent number} of $v$?
\item
In a recent paper by Bj\"orner and Ekedahl \cite{bjornerekedahl}
it is shown (for any crystallographic Coxeter group)
that $0\le i<j\le\ell(w)-i$ implies
$f_i^w\le f_j^w$, where $f_i^w$ is the $q^i$-coefficient
of $\poin_{[\id,w]}(q)$. Perhaps one can say more about
the particular Poincar\'e polynomials discussed in the present paper.
Are they unimodal, for instance?
\item
As noted by Gasharov and Reiner \cite[p.~559]{gasharovreiner},
Ding's partition varieties~\cite{ding} correspond
not only to certain Bruhat intervals of the whole group
$A_n$ but also to some intervals of the quotient $A_n/(A_n)_J$ for
certain parabolic subgroups $(A_n)_J$.
Can something similar be done in the more general setting of
skew partition varieties?
\item
Given a polynomial, what board (i.e.~zero-one matrix), if any,
has it as its rook polynomial?
In a recent paper \cite{mitchell} Mitchell showed that distinct
increasing Ferrers boards have distinct rook polynomials.
What is true for skew Ferrers boards?
\item
Develin~\cite{develin} classified the isomorphism types of
lower Bruhat intervals of 312-avoiding permutations by using
the connection to rook posets discovered by Ding~\cite{ding}.
What are the isomorphism types of lower Bruhat intervals
of permutations that avoid the patterns
4231, 35142, 42513, and 351624?
\end{itemize}

\section{Acknowledgement}
I would like to thank Torsten Ekedahl for asking a question that
made me write this paper. Many thanks also to Anders Bj\"orner
for introducing me to Coxeter groups and Schubert varieties.

\bibliographystyle{habbrv}
\bibliography{bruhat.bib}

\end{document}